\documentclass[a4paper,10pt]{amsart}

\usepackage{amsmath,amsthm,amssymb,amscd}
\usepackage[latin1]{inputenc}
\usepackage[T1]{fontenc}
\usepackage[english]{babel}
\usepackage[all]{xy}

\makeatletter
 \newtheorem*{def*}{Definition}
  
 \newtheorem*{thm*}{Theorem}
 \newtheorem{prop}{Proposition}

 \theoremstyle{definition}
 
 \newtheorem*{exp*}{Example}
 
 \newtheorem*{rem*}{Remark}

\title{Universal Prolongation of Linear Partial Differential Equations on Filtered Manifolds}

\author{Katharina Neusser}
\address{Katharina Neusser, Faculty of Mathematics, University of Vienna\\ Nordbergstra\ss e 15\\ A-1090 Wien\\
Austria} 
\email{katharina.neusser@univie.ac.at}

\thanks{This work was supported by Initiativkolleg IK-1008 of the University of Vienna and KWA grant of the University of Vienna}

\begin{document}

\maketitle

\begin{abstract}
The aim of this article is to show that systems of linear partial differential equations on filtered manifolds, which are of weighted finite type, can be canonically rewritten as first order systems of a certain type.
\\This leads immediately to obstructions to the existence of solutions. Moreover, we will deduce that the solution space of such equations is always finite dimensional.
\end{abstract}

\section{Introduction}
A filtered manifold is  a smooth manifold $M$ together with a filtration of the tangent bundle $TM=T^{-k}M\supset...\supset T^{-1}M$ by smooth subbundles such that the Lie bracket $[\xi,\eta]$ of a section
$\xi$ of $T^{i}M$ and a section $\eta$ of $T^jM$ is a section of $T^{i+j}M$. 
To each point $x\in M$ one can associate the graded vector space $\textrm{gr}(T_xM)=\bigoplus T^{i}_xM/T^{i+1}_xM$. The Lie bracket of vector fields induces a Lie bracket on this graded vector space, which makes 
$\textrm{gr}(T_xM)$ into a nilpotent graded Lie algebra. This graded nilpotent Lie algebra should be seen as the linear first order approximation to the filtered manifold at the point $x$.
\\Studying differential equations on filtered manifolds it turns out that, in addition to replacing the usual tangent space at $x$ by the graded nilpotent Lie algebra $\textrm{gr}(T_xM)$, one should also change the notion of order of differential operators according to the filtration of the tangent bundle. 
One can view a contact structure $TM=T^{-2}M\supset T^{-1}M$ as a filtered manifold structure. In this special case, this means that a derivative in direction transversal to the contact subbundle $T^{-1}M$ should be considered as an operator of order two rather than one. This leads to a notion of symbol for  differential operators on $M$ which fits naturally together with the contact structure and which can be considered as the principal part of such operators. In the context of contact geometry the idea to study differential operators on $M$ by means of their \textit{weighted} symbol goes back to the 70's and 80's of the last century and is usually referred to as Heisenberg calculus, cf. \cite{BG} and \cite{Tay}.
\\Independently of these developments in contact geometry, T. Morimoto started in the 90's to study differential equations on general filtered manifolds and developed a formal theory, cf. \cite{Morimoto1}, \cite{Morimoto2} and \cite{Morimoto3}. By adjusting the notion of order of differentiation to the filtration of a filtered manifold, he introduced a concept of weighted jet bundles, which provides the convenient framework to study differential operators between sections of vector bundles over a filtered manifold.
This leads to a notion of symbol which can be naturally viewed as the principal part of differential operators on filtered manifolds.
In \cite{Morimoto1} and \cite{Morimoto2} Morimoto also established an existence theorem for analytic solutions for certain differential equations on filtered manifolds.
\\In \cite{Spencer} Spencer extracts an important class of systems of differential equations, namely those, for which after a finite number of prolongations their symbol vanishes. He calls these systems of finite type. Such equations have always a finite dimensional solution space, since a solution is already determined by a finite jet in a single point. 
Working in the setting of weighted jet bundles one can analogously define systems of differential equations of weighted finite type. If our manifold is trivially filtered, the weighted jet bundles are the usual ones and the definition of weighted finite type coincides with the notion of Spencer.  
\\In this article we will study systems of linear differential equations on filtered manifolds by looking at their weighted symbols. Our aim is to show that to a system of weighted finite type one can always canonically associate a differential operator of weighted order one with injective weighted symbol whose kernel describes the solutions of this system. From this one can immediately see that also in the case of a differential equation of weighted finite type the solution space is always finite dimensional. Additionally, we will see that rewriting our equation in this way as a system of weighted order one leads directly to algebraic obstructions to the existence of solutions.
\\\\\textbf{Acknowledgments}
\\This article evolved from discussions with Micheal Eastwood, who advised me of some aspects of the work of Hubert Goldschmidt and Donald C. Spencer on differential equations. 
\\I would also like to thank Andreas \v Cap and Tohru Morimoto for helpful discussions.

\section{Weighted Jet Spaces}
In this section we recall some basic facts about filtered manifolds and discuss the concept of weighted jet bundles of sections
of vector bundles over filtered manifolds as it was introduced by Morimoto in \cite{Morimoto1} in order to study differential equations on filtered manifolds. For a more detailed discussion about differential equations on filtered manifolds and weighted jet bundles see also \cite{Morimoto3}.
\subsection{Filtered Manifolds}
As already mentioned, by a filtered manifold we understand a smooth manifold $M$ together with a filtration $TM=T^{-k}M\supset...\supset T^{-1}M$ of the tangent bundle by smooth subbundles, which is compatible with the Lie bracket of vector fields. 
Compatibility with the Lie bracket of vector fields here means that for sections $\xi$ of $T^iM$ and $\eta$ of $T^jM$ the Lie bracket $[\xi,\eta]$ is a section of $T^{i+j}M$, where we set $T^jM=TM$ for $j\leq-k$ and $T^jM=0$ for $j\geq0$.
\\Given a filtered manifold $M$ one can form the associated graded vector bundle $\textrm{gr}(TM)$. It is obtained 
by taking the pairwise quotients of the filtration components of the tangent bundle $$\textrm{gr}(TM)=\bigoplus_{i={-k}}^{-1}T^iM/T^{i+1}M.$$ We set $\textrm{gr}_{-i}(TM)=T^{-i}M/T^{-i+1}M$.
\\Now consider the operator $\Gamma(T^iM)\times\Gamma(T^jM)\rightarrow \textrm{gr}_{i+j}(TM)$ given by $(\xi,\eta)\mapsto q([\xi,\eta])$, where $q$ is the projection from $T^{i+j}M$ to $\textrm{gr}_{i+j}(TM)$. 
This operator is bilinear over smooth functions and therefore induced by a bundle map $T^iM\times T^jM\rightarrow \textrm{gr}_{i+j}(TM)$. Moreover it obviously factorizes to a bundle map $\textrm{gr}_i(TM)\times\textrm{gr}_j(TM)\rightarrow\textrm{gr}_{i+j}(TM)$, since for $\xi\in\Gamma(T^{i+1}M)$ we have $[\xi,\eta]\in\Gamma(T^{i+j+1}M)$.
Hence we obtain a tensorial bracket $$\{ \quad,\quad \}:\textrm{gr}(TM)\times\textrm{gr}(TM)\rightarrow\textrm{gr}(TM)$$ 
on the associated graded bundle which makes each fiber $\textrm{gr}(T_xM)$ over some point $x$ into a nilpotent graded Lie algebra.
The nilpotent graded Lie algebra $\textrm{gr}(T_xM)$ is called the symbol algebra of the filtered manifold at the point $x$. 
\\Suppose $f:M\rightarrow M$ is a local diffeomorphism whose tangent map preserves the filtration. Then the tangent map at each point $x\in M$ induces a linear isomorphism between $\textrm{gr}(T_xM)$ and $\textrm{gr}(T_{f(x)}M)$, and compatibility of $f$ with the Lie bracket easily implies that this actually is an isomorphism of Lie algebras. Hence the symbol algebra at $x$ should be seen as the first order linear approximation to the manifold at the point $x$. In general the symbol algebra may change from point to point. However, we will always assume that $\textrm{gr}(TM)$ is locally trivial as a bundle of Lie algebras. 
\\In addition, we will assume that $\textrm{gr}(T_xM)$ is generated as Lie algebra by $\textrm{gr}_{-1}(T_xM)$. 
This means that the whole filtration is determined by the subbundle $T^{-1}M$ and the filtration is just a neat way to encode the non-integrability properties of this subbundle. 
\\In the case of a trivial filtered manifold $TM=T^{-1}M$ the associated graded is just the tangent bundle, where the tangent space at each point is viewed as an abelian Lie algebra.
\begin{exp*}
Suppose $M$ is a smooth manifold of dimension $2n+1$ endowed with a contact structure, i.e. a maximally non-integrable distribution $T^{-1}M\subset TM$ of rank $2n$.
The fibers $\textrm{gr}_x(TM)=T^{-1}_xM\oplus T_xM/T^{-1}_xM$  of the associated graded bundle are then isomorphic to the Heisenberg Lie algebra $\mathfrak h=\mathbb R^{2n}\oplus \mathbb R$.  
\end{exp*}
\subsection{Differential operators and weighted jet bundles}
Studying analytic properties of differential operators on some manifold $M$, one may first look at the symbols of these operators.
If $M$ is a filtered manifold, it turns out that the usual symbol is not the appropriate object to consider and it should be replaced by a notion of symbol that reflects the geometric structure on $M$ given by the filtration on the tangent bundle. 
\\Let us consider an example. Suppose $M$ is the Heisenberg group $\mathbb R^{2n+1}$ endowed with its canonical contact structure $T^{-1}M\subset TM$. 
Denoting by $(x_1,..,x_n,y_1,..,y_n,z)$ the coordinates on $\mathbb R^{2n+1}$ and by $X_1,...,X_n,Y_1,...,Y_n,Z$ the right invariant vector fields, the distribution is spanned in each point by $X_i, Y_i$ for $i=1,...,n$.
Now consider the following operator acting on smooth functions on $M$ $$D=\sum_{j=1}^n-X_j^2-Y^2_j+iaZ \textrm{ with } a\in\mathbb C.$$ 
It can be shown, cf. \cite{VanErp}, that the analytic properties of this operator highly depend on the constant $a$. However, this can never be read off from the usual symbol, since the term $iaZ$ is not part of it. This suggests that a derivative transversal to the contact distribution should have rather order two than one to obtain a notion of symbol that includes the term $iaZ$. 
\\In the case of a general filtered manifold the situation is similar. Once one has replaced the role of the usual tangent space at some point $x\in M$ by the symbol algebra at that point, one should also adjust the notion of order of differentiation according to the filtration of the tangent bundle, in order to obtain a notion of symbol that can be seen as representing the principal parts of operators on $M$. 
\begin{def*}
(1) A local vector field $\xi$ of $M$ is of weighted order $\leq r$, if $\xi\in\Gamma(T^{-r}M)$. 
The minimum of all such $r$ is called the weighted order $\textrm{ord}(\xi)$ of $\xi$.
\\(2) A linear differential operator $D:C^{\infty}(M,\mathbb R)\rightarrow C^{\infty}(M,\mathbb R)$ between smooth functions from $M$ to $\mathbb R$
is of weighted order $\leq r$, if it can be locally written as 
$$D=\sum_i \xi_{i_1}...\xi_{i_{s(i)}}$$ for local vector fields
$\xi_{i_j}$ with 
$\sum_{\ell=1}^{s(i)} \textrm{ord}(\xi_{i_\ell})\leq r$ for all $i$.  
\\The minimum of all such $r$ is then the weighted order of $D$.  
 \end{def*}
Suppose $E\rightarrow M$ is a smooth vector bundle of constant rank over a filtered manifold $M$ and denote by $\Gamma_x(E)$ the space of germs of sections of $E$ at the point $x\in M$.
\\Then we define two sections $s,s'\in\Gamma_x(E)$ to be $r$-equivalent $\sim_r$
 if  $$D(\langle\lambda,s-s'\rangle)(x)=0$$
for all differential operators $D$ on $M$ of weighted order $\leq r$ and all sections $\lambda$ of the dual bundle $E^*$, where 
$\langle$  ,  $\rangle: \Gamma(E^*)\times\Gamma(E)\rightarrow C^\infty(M,\mathbb R)$ is the evaluation.
\\The space of weighted jets of order $r$ with source $x\in M$ is then defined as the quotient space $$J_x^r(E):=\Gamma_x(E)/\sim_r.$$
\\For $s\in\Gamma_x(E)$ we denote by $j^r_xs$ the class of $s$ in $J^r_x(E)$.
\\The space of weighted jets of order $r$ is given by taking the disjoint union over $x$ of $J^r_x(E)$  $$J^r(E):=\bigsqcup_{x\in M}J_x^r(E)$$
and we have a natural projection $J^r(E)\rightarrow M$. 
\\It is not difficult to see that for every vector bundle chart of $E$ one can construct a local trivialization of $J^r(E)$. Hence we can endow $J^r(E)$ with the unique  manifold structure such that $J^r(E)\rightarrow M$ is a vector bundle and these trivializations become smooth vector bundle charts.
\\The natural projections $$\pi^r_s:J^r(E)\rightarrow J^s(E) \textrm{ for } r>s$$ are then easily seen to be vector bundle homomorphisms.
\\Denoting by $j^r:\Gamma(E)\rightarrow\Gamma(J^r(E))$ the universal differential operator of weighted order $r$ given by 
$s\mapsto(x\mapsto j^r_xs)$, we can define the weighted order of a general linear differential operator as follows.

\begin{def*}
Suppose $E$ and $F$ are vector bundles of constant rank over a filtered manifold $M$.
A linear differential operator $D:\Gamma(E)\rightarrow\Gamma(F)$ between sections of $E$ and $F$ is of weighted order $\leq r$, if there exists a vector bundle map $\phi: J^r(E)\rightarrow F$ such that $D=\phi\circ j^r$.
The smallest integer $r$ such that this holds, is called the weighted order of $D$.
\end{def*}

Since, of course, every bundle map $\phi:J^r(E)\rightarrow F$ defines a differential operator of order $r$ by $D=\phi\circ j^r$, we can equivalently view a differential operator of weighted order $r$ as a bundle map from $J^r(E)$ to $F$.
\subsection{Weighted symbols of differential operators}
The symbol $\sigma(\phi)$ of a differential operator $\phi: J^{r}(E)\rightarrow F$ is the restriction of $\phi$ to the kernel of the projection $\pi^r_{r-1}: J^r(E)\rightarrow J^{r-1}(E)$. Let us describe this kernel more explicitly:
\\The associated graded bundle $\textrm{gr}(TM)$ is a vector bundle of nilpotent graded Lie algebras. So one can consider the universal enveloping algebra $\mathcal{U}(\textrm{gr}(T_xM))$ of the Lie algebra $\textrm{gr}(T_xM)$ defined by 
$$\mathcal{U}(\textrm{gr}(T_xM))=\mathcal T(\textrm{gr}(T_xM))/I$$
where $\mathcal T(\textrm{gr}(T_xM))$ is the tensor algebra of $\textrm{gr}(T_xM)$ and $I$ is the ideal generated by elements of the form $X\otimes Y-Y\otimes X-\{X,Y\}$ for $X, Y\in\textrm{gr}(T_xM)$.
\\The grading of $\textrm{gr}(T_xM)$ induces an algebra grading on the tensor algebra $\mathcal T(\textrm{gr}(T_xM))$ as follows: 
An element $X_1\otimes...\otimes X_\ell\in \mathcal T(\textrm{gr}(T_xM))$ is defined to be of degree $s$, if  $\sum_i \textrm{deg}(X_i)=s$ with $\textrm{deg}(X_i)=p$ for $X_i\in\textrm{gr}_p(T_xM)$. Since the ideal $I$ is homogeneous, this grading factorizes to an algebra grading on the universal enveloping algebra
$$\mathcal{U}(\textrm{gr}(T_xM))=\bigoplus_{i\leq 0}\mathcal{U}_i(\textrm{gr}(T_xM)).$$
The disjoint union $\bigsqcup_{x\in M}\mathcal{U}_{i}(\textrm{gr}(T_xM))$ is easily seen to be a vector bundle over $M$, which we denote by $\mathcal{U}_{i}(\textrm{gr}(TM))$.
\\Let us now consider the kernel of the projection $\pi^r_{r-1}:J^r(E)\rightarrow J^{r-1}(E)$.
Suppose $s$ is a local section with 
$j_x^{r-1}s=0$ for some point $x$ and take some local trivialization to write $s$ as $(s_1,...,s_n):U\subseteq M\rightarrow \mathbb R^n$ ($\textrm{rank}(E)=n$) with $x\in U$. 
For $\xi_1,..,\xi_\ell\in T_xM$ with $\sum_i \textrm{ord}(\xi_i)=r$ we have
the multilinear map $$(\xi_1,...,\xi_\ell)\mapsto\xi_1\cdot...\cdot\xi_\ell\cdot s$$ from $T_xM\times...\times T_xM$ to $\mathbb R^n$ given by iterated differentiation.
Since the $r-1$ jet of $s$ at $x$ vanishes, this map induces a linear map $\mathcal T_{-r}(\textrm{gr}(T_xM))\rightarrow\mathbb R^n$.
Additionally we have the symmetries of differentiation, like for example $\xi_1\cdot\xi_2\cdot...\cdot\xi_\ell\cdot s-\xi_2\cdot\xi_1\cdot...\cdot\xi_\ell\cdot s=[\xi_1,\xi_2]\cdot...\cdot\xi_\ell\cdot s$, which equals $\{\xi_1,\xi_2\}\cdot...\cdot\xi_\ell\cdot s$, since $j^{r-1}_xs=0$. Hence $\mathcal T_{-r}(\textrm{gr}(T_xM))\rightarrow\mathbb R^n$ factorizes to a linear map $\mathcal{U}_{-r}(\textrm{gr}(T_xM))\rightarrow\mathbb R^n$.
Via the chosen trivialization, any element of the kernel of $\pi^r_{r-1}$ determines an element in $\mathcal{U}_{-r}(\textrm{gr}(TM))^*\otimes E$ and vice versa. It is easy to see that this is independent of the chosen trivialization and
so we get the natural exact sequence of vector bundles: 
\[
\begin{CD}
0 @>>> \mathcal{U}_{-r}(\textrm{gr}(TM))^*\otimes E @>\iota>> J^r(E) @>\pi^r_{r-1}>> J^{r-1}(E)@>>> 0
\end{CD}
\]

Hence the symbol of an operator $\phi:J^r(E)\rightarrow F$ can be viewed as a bundle map $$\sigma({\phi}):\mathcal{U}_{-r}(\textrm{gr}(TM))^*\otimes E\rightarrow F.$$
\begin{rem*}
If $M$ is a trivial filtered manifold $TM=T^{-1}M$, the bundle $J^r(E)$ is just the usual bundle of jets of order $r$ and we obtain the usual notion of symbol for differential operators. 
The universal enveloping algebra of the abelian algebra $\textrm{gr}(T_xM)=T_xM$ coincides with the symmetric algebra of $T_xM$.
\end{rem*}

\section{Universal Prolongation of Linear Differential Equations on filtered manifolds}
A differential operator $\phi: J^r(E)\rightarrow F$ induces the following maps:
\\\\The $\ell$-th-prolongation $p_{\ell}(\phi):J^{r+\ell}(E)\rightarrow J^{\ell}(F)$ of $\phi$ given by
$$p_{\ell}(\phi)(j^{r+\ell}_xs)=j^{\ell}_x(\phi(j^rs)).$$
This is well defined, since the righthand side just depends on the weighted
$r+\ell$ jet of $s$ at the point $x$.
This map can be characterized as the unique vector bundle map such that the diagram
\begin{equation*}
\xymatrix{\Gamma({J^{r+\ell}(E)})\ar[r]^{{p_{\ell}(\phi)}}& \Gamma({J^{\ell}(F)})\\\Gamma({E})\ar[r]^{{\phi\circ j^r}}\ar[u]^{j^{r+\ell}}& \Gamma({F})\ar[u]^{j^{\ell}}\\}
\end{equation*}
commutes .
\\In particular, we have the bundle map $$p_{\ell}(id_r): J^{r+\ell}(E)\rightarrow J^{\ell}(J^r(E))$$ where $id_r$ is the identity map on $J^r(E)$. 
Any derivative in direction transversal to $T^{-1}M$ can be expressed by iterated derivatives in direction of the the subbundle $T^{-1}M$, since we assumed that $\textrm{gr}_{-1}(T_xM)$ generates $\textrm{gr}(T_xM)$ as Lie algebra. Therefore this vector bundle map is injective.
\\The operator $\phi$ induces also a vector bundle map $e_\ell(\phi): J^\ell(J^r(E))\rightarrow J^\ell(F)$
defined by $$e_\ell(j^{\ell}_xs)=j^{\ell}_x(\phi(s)).$$
It is the unique vector bundle map such that the diagram 
\begin{equation*}
\xymatrix{\Gamma(J^{\ell}(J^r(E)))\ar[r]^{{e_{\ell}(\phi)}}& \Gamma(J^{\ell}(F))\\\Gamma(J^r(E))\ar[r]^{{\phi}}\ar[u]^{j^\ell}& \Gamma(F)\ar[u]^{j^\ell}\\}
\end{equation*}
commutes.
By definition we have $e_\ell(\phi)\circ p_\ell(id_r)=p_\ell(\phi)$.
\\\\Since we have the inclusion $p_1(id_r):J^{r+1}(E)\hookrightarrow J^1(J^{r}(E))$, we can consider the operator $\delta^r$ of weighted order one defined
by the projection $$J^1(J^r(E))\rightarrow J^1(J^r(E))/J^{r+1}(E).$$
This operator can now be characterized, analogously as in \cite{Gold2} for usual jet bundles:
\\We have the following commutative exact diagram
\begin{equation*}
\xymatrix{
   & 0\ar[d] & 0\ar[d]& & \\0\ar[r]^{}&\mathcal{U}_{-r-1}(\textrm{gr}(TM))^*\otimes E\ar[r]\ar[d]_{\iota}&\textrm{gr}_{-1}(TM)^*\otimes J^r(E)\ar[r]^{}\ar[d]_{\iota}&W^r\ar[r] &0\\0 \ar[r]& J^{r+1}(E)\ar[r]^{p_1(id_r)}\ar[d]_{\pi^{r+1}_r}& J^1(J^r(E))\ar[r]^{\delta^r}\ar[d]_{\pi^1_0}& J^1(J^r(E))/J^{r+1}(E)\ar[r]\ar[d]&0\\0\ar[r]& J^r(E)\ar[r]\ar[d]&J^r(E)\ar[r]\ar[d]& 0 & \\ & 0 & 0 & & \\ }
\end{equation*}
where the inclusion of $\mathcal{U}_{-r-1}(\textrm{gr}(TM))^*\otimes E$ into $\textrm{gr}_{-1}(TM)^*\otimes J^r(E)$ is obtained by the commutativity of the next two rows and the space $W^r$ is defined by the diagram.
Moreover, this diagram induces an isomorphism of vector bundles between $W^r$ and $J^1(J^r(E))/J^{r+1}(E)$.
Therefore we can view $\delta^r$ as an operator from $J^1(J^r(E))$ to $W^r$.
Hence we have the following proposition: 
\begin{prop}
There exists a unique differential operator $$\delta^r: J^1(J^r(E))\rightarrow W^r$$ of weighted order one
such that
\\$\bullet$ the kernel of $\delta^r$ is $J^{r+1}(E)$ 
\\$\bullet$ the symbol $\sigma(\delta^r):\textrm{\emph{gr}}_{-1}(TM)^*\otimes J^r(E)\rightarrow W^r$ is the projection.
\end{prop}
\begin{proof}
The uniqueness follows from the exactness of the diagram above.
\end{proof}
By a similar reasoning as in \cite{Gold2} we can now deduce the existence of a first order operator 
$S: J^1(J^r(E))\rightarrow \textrm{gr}_{-1}(TM)^*\otimes J^{r-1}(E)$. We will call $S$ the weighted Spencer operator.
\begin{prop}
There exists a unique differential operator $$S:J^1(J^r(E))\rightarrow\textrm{\emph{gr}}_{-1}(TM)^*\otimes J^{r-1}(E)$$ of weighted order one such that
\\$\bullet$ $J^{r+1}(E)\subseteq \ker(S)$ 
\\$\bullet$ the symbol $\sigma(S):\textrm{\emph{gr}}_{-1}(TM)^*\otimes J^{r}(E)\rightarrow\textrm{\emph{gr}}_{-1}(TM)^*\otimes J^{r-1}(E)$ is $id\otimes\pi^r_{r-1}$.
\\Moreover, we have the following exact sequence of sheaves:
\[
\begin{CD}
0 @>>> \Gamma(E) @>j^r>> \Gamma(J^r(E)) @>S\circ j^1>> \Gamma(\textrm{\emph {gr}}_{-1}(TM)^*\otimes J^{r-1}(E)) @>>> 0
\end{CD}
\] 
\end{prop}
\begin{proof}
If such an operator exists, it must factorize over $\delta^r$ by proposition $1$, since $J^{r+1}(E)\subseteq \ker(S)$. This means that it
has to be of the form $S=\psi\circ\delta^r$ for some bundle map $\psi:W^r\rightarrow\textrm{gr}_{-1}(TM)^*\otimes J^{r-1}(E)$. By the second property $\psi$ has to satisfy that $\sigma(S)=\psi\circ \sigma(\delta^r)$ equals the projection $id\otimes\pi^r_{r-1}$.
To see that such a map $\psi$ exists, we have to show that $id\otimes\pi^{r}_{r-1}$ factorizes over $\sigma(\delta^r)$.
\\We already know that $\ker(\sigma(\delta^r))=\ker(\pi^{r+1}_r)$, which is mapped under the inclusion $p_1(id_r): J^{r+1}(E)\hookrightarrow J^1(J^r(E))$ to $\textrm{gr}_{-1}(TM)^*\otimes J^r(E)$. 
Since the map $e_1(\pi^r_{r-1}): J^1(J^r(E))\rightarrow J^1(J^{r-1}(E))$ has symbol $\iota\circ id\otimes\pi^r_{r-1}$ and we have the following commutative diagram
\[
\begin{CD}
J^{r+1}(E) @>p_1(id_r)>> J^1(J^r(E))\\
@V\pi^{r+1}_rVV  @VVe_1(\pi^r_{r-1})V\\
J^r(E) @>p_1(id_{r-1})>> J^1(J^{r-1}(E))
\end{CD}
\] we conclude that $\ker(\pi^{r+1}_r)$ is mapped under the inclusion $p_1(id_r)$ to the kernel of $id\otimes\pi^r_{r-1}$. 
Hence $id\otimes\pi^r_{r-1}$ factorizes over $\sigma(\delta^r)$. So there exists a unique bundle map $\psi: W^r\rightarrow\textrm{gr}_{-1}(TM)^*\otimes J^{r-1}(E)$ with $\psi\circ \sigma(\delta^r)=id\otimes\pi^r_{r-1}$ and we can define $S=\psi\circ\delta^r$.
\\To show the exactness of the sequence above let us describe $S$ in another way. Consider the bundle map $$e_1(\pi^r_{r-1})-p_1(id_{r-1})\circ\pi^1_0:J^1(J^r(E))\rightarrow J^1(J^{r-1}(E)).$$
Since $\pi^1_0\circ e_1(\pi^r_{r-1})=\pi^1_0\circ p_1(id_{r-1})\circ\pi^1_0$, this operator actually has values in 
$\textrm{gr}_{-1}(TM)^*\otimes J^{r-1}(E)$. Moreover, $J^{r+1}(E)$ lies in its kernel by the commutative diagram above and the symbol is given by the symbol of $e_1(\pi^r_{r-1})$ which equals $\iota\circ id\otimes\pi^r_{r-1}$.
Hence viewing $S$ as an operator from $J^1(J^r(E))$ to $J^1(J^{r-1}(E))$ by means of the inclusion $\iota:\textrm{gr}_{-1}(TM)^*\otimes J^{r-1}(E)\hookrightarrow J^1(J^{r-1}(E))$, we must have $$S=e_1(\pi^r_{r-1})-p_1(id_{r-1})\circ\pi^1_0:J^1(J^r(E))\rightarrow J^1(J^{r-1}(E)).$$
Suppose now we have a section of $J^r(E)$ which can be written as $j^rs$ for some $s\in\Gamma(E)$. Then it lies in the kernel of 
$S\circ j^1$, since $e_1(\pi^r_{r-1})\circ j^1(j^rs)=j^1(\pi^r_{r-1}(j^rs))=j^1(j^{r-1}s)=p_1(id_{r-1})(j^{r}s)$.
\\To show the converse one can proceed by induction on $r$. 
\\If $r=1$, then for $s\in\Gamma(J^1(E))$ to be in the kernel of $S\circ j^1$ means $j^1(\pi^1_0(s))=p_1(id_0)s=s$. 
Now suppose the assertion holds for $r$. If $s\in\Gamma(J^{r+1}(E))$ satisfies $j^1(\pi^{r+1}_{r}(s))=p_1(id_r)(s)$, then $e_1(\pi^r_{r-1})(j^1(\pi^{r+1}_{r}(s)))=e_1(\pi^r_{r-1})(p_1(id_r)(s))$. From the commutative diagram above we know that the right side coincides with $p_1(id_{r-1})(\pi^{r+1}_r(s))$. By the induction hypothesis $\pi^{r+1}_r(s)=j^r(u)$ for some $u\in\Gamma(E)$. Now $s$ must equal $j^{r+1}u$, since $p_1(id_r)(s)=j^1(\pi^{r+1}_{r}(s))=j^1(j^r u)$
and $p_1(id_r)$ is injective.
\end{proof}
Now we want to study systems of linear differential equations on filtered manifolds.
Suppose we have vector bundle map 
$$\phi:J^r(E)\rightarrow F$$
of constant rank, then the subbundle of $J^r(E)$ defined by its kernel
$$Q^r:=\ker(\phi)$$ is called the linear system of differential equations associated to the operator $\phi$. 
A solution of $Q^r$ is a section $s$ of $E$ satisfying $\phi(j^rs)=0$.
\\\\The $\ell$-th prolongation $Q^{r+\ell}$ of $Q^r$ is defined to be the kernel of $p_\ell(\phi):J^{r+\ell}(E)\rightarrow J^\ell(F)$.
Since the diagram
\begin{equation*}
\xymatrix{ & & 0\ar[d] & 0\ar[d]\\ & & J^{r+\ell}(E)\ar[r]^{p_\ell(\phi)}\ar[d]_{p_\ell(id_r)}& J^\ell(F)\ar[d]^{id}\\0\ar[r]& J^\ell(Q^r)\ar[r]& J^\ell(J^r(E))\ar[r]^{e_\ell(\phi)}& J^\ell(F)\\}
\end{equation*}
commutes, we have $$Q^{r+\ell}=J^\ell(Q^r)\cap J^{r+\ell}(E).$$
\\In general, the bundle map $p_\ell(\phi)$ is not of constant rank and $Q^{r+\ell}$ need not to be a vector bundle.
We call $\phi:J^r(E)\rightarrow F$ regular, if $p_{\ell}(\phi)$ is of constant rank for all $\ell\geq 0$.
\\The symbol of the prolonged equation $Q^{r+\ell}$ is the family of vector spaces $g^{r+\ell}:=\{g^{r+\ell}_x\}_{x\in M}$ over $M$, where $g^{r+\ell}_x$ is the kernel of the linear map $Q^{r+\ell}_x\rightarrow Q^{r+\ell-1}_x$ given by the restriction of the projection $\pi^{r+\ell}_{r+\ell-1}$ to $Q^{r+\ell}_x$.
\\For all $\ell\geq 1$ we have a bundle map $\sigma_\ell(\phi):\mathcal{U}_{-r-\ell}(\textrm{gr}(TM))^*\otimes E\rightarrow\mathcal{U}_{-\ell}(\textrm{gr}(TM))^*\otimes F$, which we call the $\ell$-th symbol mapping. It is defined by the following (fiberwise) commutative diagram:
\begin{equation*}
\xymatrix{ & 0\ar[d]& 0\ar[d] & 0\ar[d]\\ 0\ar[r]& g_{r+\ell}\ar[r]\ar[d]&\mathcal U_{-r-\ell}(\textrm{gr}(TM))^*\otimes E\ar[r]^{\sigma_\ell(\phi)}\ar[d]&\mathcal U_{-\ell}(\textrm{gr}(TM))^*\otimes F\ar[d]\\0\ar[r]& Q^{r+\ell}\ar[r]\ar[d]& J^{r+\ell}(E)\ar[r]^{p_\ell(\phi)}\ar[d]& J^\ell(F)\ar[d]\\0\ar[r]& Q^{r+\ell-1}\ar[r]& J^{r+\ell-1}(E)\ar[r]^{p_{\ell-1}(\phi)}& J^{\ell-1}(F)\\}
\end{equation*}
By definition the kernel of $\sigma_{\ell}(\phi)$ is $g^{r+\ell}$ viewed as a subset of 
$\mathcal U_{-r-\ell}(\textrm{gr}(TM))^*\otimes E$.
\\Since $Q^{r+\ell}=J^\ell(Q^r)\cap J^{r+\ell}(E)$ and the diagram
\[
\begin{CD}
J^{r+\ell}(E) @>p_\ell(id_r)>> J^\ell(J^r(E))\\
@V\pi^{r+\ell}_{r+\ell-1}VV  @VV\pi^\ell_{\ell-1}V\\
J^{r+\ell-1}(E) @>p_{\ell-1}(id_{r})>> J^{\ell-1}(J^{r}(E))
\end{CD}
\]

commutes, we conclude that $$g^{r+\ell}_x=\mathcal U_{-r-\ell}(\textrm{gr}(T_xM))^*\otimes E_x\cap\mathcal U_{-\ell}(\textrm{gr}(T_xM))^*\otimes K_x$$ where $K=\{K_x\}_{x\in M}$ is the kernel of the symbol $\sigma(\phi)$ of $\phi$.
\\\\A differential equation $Q^r\subset J^r(E)$ is called of finite type, if there exists $m\in\mathbb N$ such that $g^{r+\ell}_x=0$ for all $x\in M$ and $\ell\geq m$.
\\\\For equations of finite type we can prove the following theorem: 

\begin{thm*}
Suppose $M$ is a filtered manifold such that $\textrm{\emph{gr}}(TM)$ is locally trivial as a bundle of Lie algebras and $\text{\emph{gr}}_{-1}(T_xM)$ generates $\textrm{\emph{gr}}(T_xM)$ for all $x\in M$ and suppose $E$ and $F$ are vector bundles over $M$. 
\\Let $D:\Gamma(E)\rightarrow \Gamma(F)$ be a regular differential operator of weighted order $r$ defining a system of
differential equations of finite type.
Then for some $\ell_0\in\mathbb N$ there exists a differential operator $$D':\Gamma(Q^{r+\ell_0})\rightarrow\Gamma(W^{r+\ell_0}) \quad\textrm{of weighted order one with injective symbol}$$ such that
$s\mapsto j^{r+\ell_0}s$ induces a bijection:
$$\{s\in\Gamma(E):D(s)=0\}\leftrightarrow\{s'\in\Gamma(Q^{r+\ell_0}):D'(s')=0\}.$$
 \end{thm*}
 \begin{proof}
 Let us denote by $\phi: J^r(E)\rightarrow F$ the bundle map associated to $D$ and by $Q^r\subseteq J^r(E)$ the differential equation given by the kernel of $\phi$.
 For all $\ell\geq 0$ we can consider the operator $D^{r+\ell}: \Gamma(Q^{r+\ell})\rightarrow \Gamma(W^{r+\ell})$ of weighted order one given by the restriction of $\delta^{r+\ell}\circ j^1$ to $\Gamma(Q^{r+\ell})$.
 \\If $s\in\Gamma(E)$ is a solution $Ds=0$, then $j^{r+\ell}s\in\Gamma(Q^{r+\ell})$ and since $j^1(j^{r+\ell}s)$ is a section of $J^{r+\ell+1}(E)\subseteq J^1(J^{r+\ell}(E))$ we also have $D^{r+\ell}(j^{r+\ell}s)=0$.
 \\And conversely, if $s'$ is a section of $Q^{r+\ell}$ such that $D^{r+\ell}(s')=0$, then $j^1s'$ is a section of $J^{r+\ell+1}(E)$. Since $J^{r+\ell+1}(E)$ is contained in the kernel of the weighted Spencer operator $J^1(J^{r+\ell}(E))\rightarrow J^1(J^{r+\ell-1}(E))$,
 the section $s'$ equals $j^{r+\ell}s$ for some section $s\in\Gamma(E)$. Obviously $\pi^{r+\ell}_0(j^{r+\ell}s)=s$ then satisfies $Ds=0$.
 \\This shows that for all $\ell\geq0$ the map $j^{r+\ell}$ induces a bijection between solutions of $D$ and solutions of $D^{r+\ell}$.
 So it remains to prove that there exists some $\ell_0$ such that $D^{r+\ell_0}$ has injective symbol.
\\The symbol of $D^{r+\ell}$ is a bundle map $\mathcal{U}_{-1}(\textrm{gr}(TM))^*\otimes Q^{r+\ell}\rightarrow W^{r+\ell}$. 
 We know that the kernel of $\sigma(\delta^{r+\ell})$ is 
 $\mathcal{U}_{-r-\ell-1}(\textrm{gr}(TM))^*\otimes E$. 
 \\Since $D^{r+\ell}$ is just the restriction of $\delta^{r+\ell}\circ j^1$ to $\Gamma(Q^{r+\ell})$, we obtain that $$\ker(\sigma(D^{r+\ell}))_x=\mathcal U_{-r-\ell-1}(\textrm{gr}(T_xM))^*\otimes E_x\cap\mathcal U_{-1}(\textrm{gr}(T_xM))^*\otimes g^{r+\ell}_x.$$ 
 But $g^{r+l}_x=\mathcal U_{-r-l}(\textrm{gr}(T_xM))^*\otimes E_x\cap\mathcal U_{-l}(\textrm{gr}(T_xM))^*\otimes K_x$ where $K$ is the kernel of the symbol of $D$. 
 Therefore we get  $$\ker(\sigma(D^{r+\ell}))_x=\mathcal U_{-r-\ell-1}(\textrm{gr}(T_xM))^*\otimes E_x\cap\mathcal U_{-\ell-1}(\textrm{gr}(T_xM))^*\otimes K_x$$ which coincides with $g^{r+\ell+1}_x$.
 \\Since the equation $Q^r$ is of finite type, there exists $\ell_0$ such that $g^{r+\ell_0+1}=0$ and hence $D^{r+\ell_0}:\Gamma(Q^{r+\ell_0})\rightarrow \Gamma(W^{r+\ell_0})$ is a differential operator of weighted order one with injective symbol, whose solutions are in bijective correspondence with solutions of the original equation $Q^r$.
 \end{proof}
 As a consequence of this theorem, we obtain that a solution of a regular differential equation of weighted finite type is already determined by a finite jet in a single point, since a solution of $D^{r+\ell_0}(s')=0$ is determined by its value in a single point. Hence the solution space of a differential equation of weighted finite type is always finite dimensional.
 \\Moreover, since $D^{r+\ell_0}$ is of weighted order one with injective symbol, it induces a vector bundle map $\rho: Q^{r+\ell_0}\rightarrow W^{r+\ell_0}/\mathcal{U}_{-1}(\textrm{gr}(TM))^*\otimes Q^{r+\ell_0}$. 
 \\Any solution $s'$ of $D^{r+\ell_0}$ must clearly also satisfy $\rho(s')=0$, which leads to obstructions for the existence of solutions. 
\begin{rem*}
The fact that a differential equation of weighted finite type has finite dimensional solution space was (by other means) already earlier observed by Morimoto, \cite{Mor4}.
\end{rem*}


\begin{thebibliography}{XX}
\bibitem{BG} R. Beals and P.C. Greiner: \textit{Calculus on Heisenberg manifolds}. Annals of Math. Studies, vol. 119. Princeton Univ. Press, Princeton NJ. 1988.
\bibitem{VanErp} E. van Erp: \textit{The Atiyah-Singer index formula for subelliptic operators on contact manifolds. Part 1.} To appear in Annals of Math. preprint arXiv: 0804.2490.
\bibitem{Gold1} H. Goldschmidt: \textit{Existence theorems for analytic linear partial differential equations}. Ann. Math., vol.86. 1967. p. 246-270.
\bibitem{Gold2} H. Goldschmidt: \textit{Prolongations of linear partial differential equations: A conjecture of \'Elie Cartan}. Ann scient. \'Ecole Norm. Sup. 4 s\'erie. t. 1. 1968. 417-444.
\bibitem{Morimoto1} T. Morimoto: \textit{Th$\acute{e}$or$\grave{e}$me de Cartan-K$\ddot{a}$hler dans une classe de fonctions formelles Gevrey.} C. R. Acad. Sci. Paris. 311, s$\acute e$rie A. 1990. p. 443-436.
\bibitem{Morimoto2} T. Morimoto:\textit{ Th$\acute{e}$or$\grave{e}$me d'existence de solutions analytiques pour des syst$\grave{e}$mes d'$\acute{e}$quations aux d$\acute{e}$riv$\acute{e}$es partielles non-lin$\acute{e}$aires avec singularit$\acute{e}$s}. C.R. Acad. Sci. Paris. 321, s$\acute{e}$rie 1. 1995. p.1491-1496.
 \bibitem{Morimoto3} T. Morimoto: \textit{Lie algebras, geometric
  structures and differential equations on filtered manifolds}. In ``Lie
  Groups Geometric Structures and Differential Equations - One Hundred
  Years after Sophus Lie'', Adv. Stud. Pure Math., Math.
  Soc. of Japan, Tokyo. 2002. 205--252.
\bibitem{Mor4}T. Morimoto: private communication.
\bibitem{Spencer} D.C. Spencer: \textit{Overdetermined systems of linear partial differential equations}. Bull. Amer. Math. Soc. 75. 1969. 179-239.
\bibitem{Tay}M.E. Taylor: \textit{Noncommutative microlocal analysis I}. Mem. Amer. Math. Soc. 52, no. 313. 1984.
\end{thebibliography}
\end{document}